\documentclass[12pt]{article}
\usepackage{bbm}
\usepackage{amsfonts}
\usepackage{pifont}
\usepackage{hyperref}
\usepackage{mathrsfs}
\usepackage{enumitem}
\usepackage{indentfirst}
\usepackage{amsmath}
\usepackage{amssymb}
\usepackage[amsmath, thmmarks]{ntheorem}
\usepackage{cite}
\usepackage{graphicx}
\usepackage{tikz}
\newtheorem{definition}{\bf Definition}[section]

\newtheorem{theorem}{\bf Theorem}[section]
\newtheorem{remark}{\bf Remark}[section]
\newtheorem{corollary}{\bf Corollary}[section]
\newtheorem{proposition}{\bf Proposition}[section]
\newtheorem{example}{\bf Example}[section]

\def\QEDopen{{\setlength{\fboxsep}{0pt}\setlength{\fboxrule}{0.2pt}\fbox{\rule[0pt]{0pt}{1.3ex}\rule[0pt]{1.3ex}{0pt}}}} 
\def\QED{\QEDopen}
\def\proof{{\bf Proof.} }
\def\endproof{\hspace*{\fill}~\QED\par\endtrivlist\unskip}

\begin{document}
\setcounter{page}{1}

\title{{\textbf{Constructing left-continuous triangular norms on complete lattices}}\thanks {Supported by the National Natural Science
Foundation of China (No. 12471440)}}
\author{Peng He\footnote{\emph{E-mail address}: 443966297@qq.com}, Xue-ping Wang\footnote{Corresponding author. xpwang1@hotmail.com; fax: +86-28-84761502}\\
\emph{\small (School of Mathematical Sciences, Sichuan Normal University}\\ \emph{\small Chengdu 610066, Sichuan, People's Republic of China)}}
\newcommand{\pp}[2]{\frac{\partial #1}{\partial #2}}
\date{}
\maketitle

\begin{quote}
{\bf Abstract}
This article focuses on the construction of left-continuous t-norms on complete lattices. The concepts of $\mathfrak{f}$-mappings and weak $\mathfrak{f}$-mappings on complete lattices are first introduced, respectively. 
They are then applied to establish the following key results: weak $\mathfrak{f}$-mappings are used to induce left-continuous t-subnorms; $\mathfrak{f}$-mappings are used to generate left-continuous t-norms whenever the top element $1$ of the complete
lattice is a completely join-irreducible element. Finally, some necessary and sufficient conditions are provided for an operator constructed by the ordinal sum of a series of 
annihilating binary operators being a left-continuous t-norm on a complete lattice.


\end{quote}

\section{Introduction}\label{intro}

Triangular norms (t-norms for short) are associative, commutative and non-decreasing binary operators 
with the neutral element $1$ on the unit interval $[0,1]$ \cite{Klement}.
As t-norms constitute special cases of compact semigroups,
the concept of ordinal sums in the sense of Clifford \cite{Clifford} provides
a method for constructing new t-norms on $[0,1]$ from given ones. The following theorem presents
the notion of an ordinal sum of t-norms and the representation of t-norms on $[0,1]$.
\begin{theorem}[\cite{Davey, Klement}] \label{00}
Let $T_i$ be a t-norm on the subinterval $[a_i, b_i]$ for each $i\in I$, where $\{(a_i, b_i)\mid i\in I\}$ is a family of nonempty, pairwise disjoint open
subintervals of $[0,1]$ and $I$ is a totally ordered index set. Then the binary operation $T: [0,1]\rightarrow [0,1]$ defined by

\begin{equation}\label{eqq0}
T(x,y)=
\begin{cases}
a_i+(b_i-a_i)T_i(\frac{x-a_i}{b_i-a_i}, \frac{y-a_i}{b_i-a_i}) &\emph{if } x,y\in [a_i, b_i],\\
\min\{x,y\} & \emph{otherwise},
\end{cases}
\end{equation}
is a t-norm.
\end{theorem}
The t-norm defined by \eqref{eqq0} is referred to as the ordinal sum of $\{T_i\mid i\in I\}$. Furthermore, every
continuous t-norm can be represented as an ordinal sum of continuous Archimedean t-norms \cite{Alsina, Klement}.

Since a lattice with the top element $1$ and the bottom element $0$ is called a bounded lattice, it constitutes a natural generalization of the unit interval $[0,1]$. Then it is quite natural to give a definition of a t-norm on a bounded lattice as follows.
\begin{definition}[\cite{De, Alsina}]\label{De2}
	\emph{A binary operator $T: L^2\rightarrow L$ with $L$ a bounded lattice is called a t-norm if, for all $a, b, c\in L$, it satisfies \\
		(a) $T(1, a)=a$;\\
		(b) if $b\leq c$ then $T(a, b)\leq T(a, c)$;\\
		(c) $T(a, b)=T(b, a)$; \\
		(d) $T(T(a, b), c)=T(a, T(b, c))$. }
\end{definition}
From Definition \ref{De2}, the pair $(L, T)$ forms an Abelian ordered semigroup with the neutral element $1$. As a consequence, $0$ is an annihilator of $T$,
and the following property holds: \\
(e) $T(x,y)\leq x\wedge y$ for all $x,y\in L$. 

The strongest and weakest t-norms on a bounded lattice $L$ are given, respectively, by\\
$$T_{M}(x, y)=x\wedge y\mbox{ and } T_{D}(x,y)=\begin{cases}
	x\wedge y &\mbox{if } 1\in\{x, y\},\\
	0 & \mbox{otherwise}.
\end{cases}$$

It is also intuitive to define ordinal sums on bounded lattices in an
analogous manner and utilize them to investigate the construction of t-norms. Saminger \cite{Sam, Sam08} constructed an ordinal sum of t-norms on a bounded lattice
by directly extending the ordinal sum construction from the unit interval $[0, 1]$.

\begin{theorem}[\cite{Sam}]
 Let $L$ be a bounded lattice, $I$ be a totally ordered index set, $\{(a_i, b_i)\mid i\in I\}$ be a family of pairwise disjoint subintervals of $L$,
 and let each $T_i$ denote a t-norm on the interval $[a_i, b_i]$.
Then the corresponding ordinal sum $T: L^2\rightarrow L$ is given by
\begin{equation}\label{eqq1}
T(x,y)=
\begin{cases}
T_i(x,y) &\emph{if } x,y\in [a_i, b_i],\\
x\wedge y & \emph{otherwise}.
\end{cases}
\end{equation}
\end{theorem}
However, the binary operator $T$ defined by \eqref{eqq1} may fail to be a t-norm on the bounded lattice $L$. In other words, constructing a t-norm that takes the form
of an ordinal sum of such interval-based t-norms on a bounded lattice is non-trivial \cite{Sam}. Therefore, many scholars have explored the construction of t-norms on bounded lattices
by imposing constraints on the lattice elements or revising Eq. \eqref{eqq1} (see e.g. \cite{Dvo, Dvo20, Gen, Lai, El, El13, Ert, Ert20, Med, Ouy}).

It is well-known that the important role of left-continuous t-norms in several fields such as probabilistic
metric spaces or fuzzy logic is undeniable \cite{Mesiar, Du, Ma, Wang02, Wang03, Weber,Fod}. However, the structure of left-continuous t-norms on complete lattices remains
an open problem. Thus, we naturally attempt to use the ordinal sum method to investigate left-continuous t-norms on complete lattices.
Unfortunately, even if the binary operators on their respective subintervals are already left-continuous t-norms,
we may still be unable to construct a left-continuous t-norm on the complete lattice---as such a t-norm may not exist as is shown by the following example.

\begin{example}\label{fl}
\emph{Consider a finite lattice $L$ depicted in Fig. 1. Note that $[0,a], [a, b]$ and $[b,1]$ are all subintervals of $L$.
It is easily verified that $T_M$ is a left-continuous t-norm on each of these subintervals.}

\emph{Suppose that there exists a left-continuous t-norm $T$ on $L$.
Then $$b=T(b,1)=T(b, c\vee e)=T(b,c)\vee T(b,e)=T(b,c)\vee 0=T(b,c).$$
Combining with $d\geq b$ and $T(d,c)\leq d\wedge c=b$, this implies that $T(d,c)=b$.
Nevertheless, we also have $$T(d, c)=T(e\vee a, c)=T(e, c)\vee T(a,c)=0\vee T(a,c)\leq a< b,$$
a contradiction. Thus, no left-continuous t-norm exists on $L$.}
\end{example}
\par\noindent\vskip50pt
\begin{minipage}{11pc}
\setlength{\unitlength}{0.7pt}\begin{picture}(600,200)
\put(290,120){\circle{6}}\put(296,120){\makebox(0,0)[l]
{\footnotesize $a$}}
\put(210,140){\circle{6}}\put(196,140){\makebox(0,0)[l]
{\footnotesize $e$}}
\put(250,80){\circle{6}}\put(246,68){\makebox(0,0)[l]
{\footnotesize $0$}}
\put(290,160){\circle{6}}\put(275,158){\makebox(0,0)[l]
{\footnotesize $b$}}
\put(290,240){\circle{6}}\put(288,250){\makebox(0,0)[l]
{\footnotesize $1$}}
\put(330,200){\circle{6}}\put(335,198){\makebox(0,0)[l]
{\footnotesize $c$}}
\put(250,200){\circle{6}}\put(235,198){\makebox(0,0)[l]
{\footnotesize $d$}}
\put(290,123){\line(0,1){34}}
\put(292,162){\line(1,1){36}}
\put(252,82){\line(1,1){36}}
\put(248,82){\line(-2,3){37}}
\put(212,142){\line(2,3){37}}
\put(288,162){\line(-1,1){36}}
\put(288,162){\line(-1,1){36}}
\put(252,202){\line(1,1){36}}
\put(328,202){\line(-1,1){36}}
\put(283,65){$L$}
\put(170,20){ Fig.1. A finite lattice $L$.}
\end{picture}
\end{minipage}

In Example \ref{fl}, the essential reason for the failure to construct a left-continuous t-norm on $L$ via binary
operators on the subintervals $[0,a], [a, b]$ and $[b,1]$ lies on the incomparability of $e$ with $a$ and $b$.
Therefore, we naturally concern whether the ordinal sum-based
construction of left-continuous t-norms is still efficient for a complete lattice in which every element is comparable with the endpoints of their respective specified subintervals. This article will positively answer the problem.

The remainder of this article is structured as follows. In Section 2, we review some fundamental concepts which will be used lately.
In Section 3, we introduce $\mathfrak{f}$-mappings and weak $\mathfrak{f}$-mappings on
complete lattices, respectively, and demonstrate that weak f-mappings can be used to induce left-continuous t-subnorms,
while f-mappings can be applied for inducing left-continuous t-norms provided that the top element $1$ of the complete
lattice is a completely join-irreducible element. In Section 4, we establish some necessary and sufficient
conditions for the operators $T$ constructed by the ordinal sum of a series of 
annihilating binary operators being left-continuous t-norms.
Conclusions are drawn in Section 5.

\section{Preliminaries}
For the sake of completeness, we review some fundamental concepts regarding lattices and t-norms on bounded lattices.\cite{Davey, Crawley73, Gratz, Birkhoff73}



An element $z$ in a lattice $L$ is said to be join-irreducible if, for all $x, y\in L$, $z=x\vee y$ implies $z=x$ or $z=y$.
For a complete lattice $L$, an element $z\in L$ is completely join-irreducible if, for every subset $X\subseteq L$, $z=\bigvee X$ implies $z\in X$.
It is straightforward that every completely join-irreducible element is join-irreducible.\cite{Crawley73}

Let $L$ be a lattice, $a, b\in L$ and $a\leq b$. The subinterval $[a, b]$ is defined as
$$[a,b]=\{x\in L\mid a\leq x\leq b \},$$
and the half-open interval $(a, b]$ is defined as
$$(a,b]=\{x\in L\mid a<x\leq b\}.$$

Let $(I,\sqsubseteq)$ be a nonempty linearly ordered index set and $\{[a_i, b_i]\mid i\in I\}$ be a family of intervals
such that for all $i, j\in I$ with $i\sqsubset j$ either $[a_i, b_i]$ and $[a_j, b_j]$ are disjoint or $b_i=a_j$. A linear sum
of $\{[a_i, b_i]\mid i\in I\}$ \cite{Davey} is the set $\bigcup_{i\in I}[a_i, b_i]$ equipped with the partial order $\leq$ defined by
$$x\leq y \Leftrightarrow (\exists i\in I, x,y\in [a_i, b_i], x\leq_i y) \mbox{ or } (\exists  i\sqsubset j\in I, x\leq b_i, y\geq a_j).$$

\begin{definition}[\cite{Mesiar}]\label{Deee2}
\emph{A binary operator $T: L^2\rightarrow L$ with $L$ a bounded lattice is called a t-subnorm if it satisfies (b), (c), (d) and (e) of Definition \ref{De2}. }
\end{definition}

Moreover, a t-subnorm is called a strong t-subnorm if $T(1,1)$=1. Obviously, a t-norm is a strong t-subnorm, but the inverse generally does not hold.

\begin{definition}[\cite{De}]\label{de002}
\emph{Let $L$ be a complete lattice and $T$ be a t-norm (resp. t-subnorm) on $L$. An element $x\in L$ is called an
idempotent element of $T$ if $T(x,x)=x$.}
\end{definition}

Let $L$ be a complete lattice. A t-norm (resp. t-subnorm) $T: L^2\rightarrow L$ is said to be left-continuous if, for all $a\in L$ and $S\subseteq L$ with $S\neq \emptyset$,
\begin{equation*}\label{Eq1}
 T(a, \vee S)=\bigvee_{s\in S} T(a,s).
\end{equation*}

In what follows, we always assume that $L$ is a complete lattice with $0$ the bottom element and $1$ the top element.
For two sets $A$ and $B$, we define $A\setminus B=\{x\in A| x\notin B\}$. In particular, if $B=\{b\}$, then we denote $A\setminus B=A\setminus b$.

\section{$\mathfrak{f}$-mappings on complete lattices}
In this section, we introduce the concepts of $\mathfrak{f}$-mappings and weak $\mathfrak{f}$-mappings on complete lattices, which are used to induce left-continuous t-norms and t-subnorms, respectively.

A complete lattice is called join-infinite distributive if, for all $a\in L$ and $S\subseteq L$ with $S\neq \emptyset$,
\begin{equation*}
	a\wedge (\bigvee S)=\bigvee_{s\in S} (a\wedge s).
\end{equation*}  
\begin{definition}\label{Def00}
\emph{Let $\mbox{id}_L$ denote an identity mapping on a complete lattice $L$. Then a mapping $f: L\rightarrow L$ is
called an $\mathfrak{f}$-mapping if, for all $x\in L$ and $S\subseteq L$ with $S\neq \emptyset$, it satisfies the following conditions:\\
(i) $f(1)=1$ (i.e., $f$ preserves the top element);\\
(ii) $f\leq \mbox{id}_L$ (i.e., $f$ is contractive, meaning $f(x)\leq x$); \\
(iii) $f^2=f$ (i.e., $f$ is idempotent, satisfying $f(f(x))=f(x)$); \\
(iv) $f(\bigvee S)=\bigvee_{s\in S} f(s)$ (i.e., $f$ preserves arbitrary joins); \\
(v) The image set of $f$, denoted by $\mbox{Im}(f)=\{f(x)\mid x\in L\}$, forms a join-infinite distributive lattice
under the partial order of $L$.}\\
\emph{A mapping $f: L\rightarrow L$ is referred to as a weak $\mathfrak{f}$-mapping if it satisfies (ii), (iii), (iv) and (v).}
\end{definition}

\begin{remark}\label{rem00}
\emph{(1) If $f$ is an $\mathfrak{f}$-mapping (resp. a weak $\mathfrak{f}$-mapping) on a complete lattice $L$, then by (iv) of Definition \ref{Def00}, $f$ is an
order-preserving mapping, i.e., $f(x)\leq f(y)$ whenever $x\leq y$.}

\emph{(2) Not every complete lattice admits an $\mathfrak{f}$-mapping. For instance, the modular lattice $M_3$ has no $\mathfrak{f}$-mapping.}

\emph{(3) If the top element $1$ of a complete lattice $L$ is a completely join-irreducible element, then there exists at least one $\mathfrak{f}$-mapping on $L$. Specifically, the mapping
\begin{equation*}
f(x)=
\begin{cases}
x &\mbox{if } x=1,\\
0 & \mbox{otherwise}
\end{cases}
\end{equation*} serves as such an $\mathfrak{f}$-mapping. Furthermore, every join-infinite distributive lattice $L$ has at least one $\mathfrak{f}$-mapping (for example, take $f=\mbox{id}_L$).}

\emph{(4) Evidently, every complete lattice $L$ has at least one weak $\mathfrak{f}$-mapping. In particular, the constant mapping $f(x)= 0$ is
the smallest weak $\mathfrak{f}$-mapping under the natural partial order.}

\emph{(5) If $f$ is an $\mathfrak{f}$-mapping (resp. a weak $\mathfrak{f}$-mapping) on a complete lattice $L$, then $\mbox{Im}(f)$ coincides with the set of fixed points of $f$ 
and $\mbox{Im}(f)$ is a join-sublattice of $L$. 
Additionally, it is worth noting that $\mbox{Im}(f)$ is not necessarily a sublattice of $L$.
To illustrate this, consider the finite lattice $L$ depicted in Fig. 2 and the mapping $f$ specified in Table \ref{Tab:00}.
It is straightforward to verify that $f$ is an $\mathfrak{f}$-mapping on $L$, however, $f(c)\wedge_L f(d)=c\wedge_L d=b\notin{\mbox{Im}}(f)$.}

\begin{table}[!h]
\centering
\caption{An $\mathfrak{f}$-mapping $f$ on $L$}
\label{Tab:00}
\begin{tabular}{c|cccccc}

 $x$ & $0$ & $a$ & $b$ & $c$ & $d$ & $1$ \\
 \hline
$f(x)$ & $0$ & $0$ & $0$ & $c$ & $d$ &$1$\\
\end{tabular}
\end{table}

\par\noindent\vskip20pt
\begin{minipage}{11pc}
\setlength{\unitlength}{0.7pt}\begin{picture}(600,250)
\put(290,120){\circle{6}}\put(296,120){\makebox(0,0)[l]
{\footnotesize $a$}}
\put(290,80){\circle{6}}\put(286,68){\makebox(0,0)[l]
{\footnotesize $0$}}
\put(290,160){\circle{6}}\put(275,158){\makebox(0,0)[l]
{\footnotesize $b$}}
\put(290,240){\circle{6}}\put(288,250){\makebox(0,0)[l]
{\footnotesize $1$}}
\put(330,200){\circle{6}}\put(335,198){\makebox(0,0)[l]
{\footnotesize $c$}}
\put(250,200){\circle{6}}\put(235,198){\makebox(0,0)[l]
{\footnotesize $d$}}
\put(290,123){\line(0,1){34}}
\put(292,162){\line(1,1){36}}
\put(290,83){\line(0,1){34}}
\put(288,162){\line(-1,1){36}}
\put(288,162){\line(-1,1){36}}
\put(252,202){\line(1,1){36}}
\put(328,202){\line(-1,1){36}}
\put(283,40){$L$}
\put(170,20){ \emph{Fig.2. A finite lattice $L$.}}
\end{picture}
\end{minipage}
\end{remark}

Let $f$ be a weak $\mathfrak{f}$-mapping on a complete lattice $L$. For convenience, in the following, 
we always use $f(x)\vee f(y)$ and $f(x)\wedge f(y)$ instead of $f(x)\vee_{\mbox{Im}(f)} f(y)$ and $f(x)\wedge_{\mbox{Im}(f)} f(y)$, respectively, where $\vee_{\mbox{Im}(f)}$ and $\wedge_{\mbox{Im}(f)}$ are the join and meet 
operations restricted to $\mbox{Im}(f)$, respectively. 

\begin{theorem}\label{pro01}
If $f$ is a weak $\mathfrak{f}$-mapping on a complete lattice $L$, then the binary operator $T: L^2\rightarrow L$ defined by $T(x, y)=f(x)\wedge f(y)$ for all $x, y\in L$
is a left-continuous t-subnorm on $L$.
\end{theorem}
\proof
The commutativity of $T$ is trivial.
For the monotonicity, let $x, y, z\in L$ with $y\leq z$. Then
$$T(x, y)=f(x)\wedge f(y)\leq f(x)\wedge f(z)=T(x, z)$$
since $f(y)\leq f(z)$, which implies that $T$ is non-decreasing.
Next, we verify the associativity. For any $x, y, z\in L$, we have
\begin{align*}
T(x, T(y,z))&=T(x, f(y)\wedge f(z))\\
&=f(x)\wedge f(f(y)\wedge f(z))\\
&=f(x)\wedge f(y)\wedge f(z) \ \ \ (\mbox{since} f^2=f \mbox{ and } \mbox{Im}(f) \mbox{ is a lattice})\\
&=f(f(x)\wedge f(y))\wedge f(z)  \ \ \ (\mbox{the same reason as above})\\
&=T(f(x)\wedge f(y), z)\\
&=T(T(x, y), z).
\end{align*}
Hence $T$ is associative. Finally, since $\mbox{Im}(f)$ is a subposet of $L$ and $f\leq \mbox{id}_L$, $T(x, y)=f(x)\wedge f(y)\leq x\wedge y$.
Combining the above results, $T$ is a t-subnorm on $L$.

It remains to verify the left-continuity. Let $x\in L$ and $\emptyset\neq S\subseteq L$. By the definition of $T$, the join-preserving property of $f$
and the join-infinite distributivity of $\mbox{Im}(f)$, we have
\begin{align*}T(x, \bigvee S)&=f(x)\wedge f(\bigvee S)\\
&=f(x)\wedge [\bigvee_{s\in S} f(s)]\\
&=\bigvee_{s\in S} (f(x)\wedge f(s))\\
&=\bigvee_{s\in S} T(x,s).
\end{align*}
Therefore, $T$ is left-continuous.

To sum up, $T$ is a left-continuous t-subnorm.
\endproof

\begin{example}\label{exx1}
\emph{Consider the lattice $L$ depicted in Fig. 2 again. Let $f$ be the mapping on $L$ defined in Table \ref{Tab:001}. Obviously, $f$ is a weak $\mathfrak{f}$-mapping.
By Theorem \ref{pro01}, we obtain a left-continuous t-subnorm $T$ on $L$, whose values are presented in Table \ref{Tab:0001}.}
\end{example}
\begin{table}[!h]
\centering
\caption{A weak $\mathfrak{f}$-mapping $f$}
\label{Tab:001}
\begin{tabular}{c|cccccc}
$x$ & $0$ & $a$ & $b$ & $c$ & $d$ & $1$ \\
 \hline
$f(x)$ & $0$ & $0$ & $b$ & $b$ & $b$ &$b$\\

\end{tabular}
\end{table}

\begin{table}[!h]
\centering
\caption{A left-continuous t-subnorm $T$}
\label{Tab:0001}
\begin{tabular}{c|cccccc}
$T$ & $0$ & $a$ & $b$ & $c$ & $d$ & $1$ \\
 \hline
$0$ & $0$ & $0$ & $0$ & $0$ & $0$ &$0$\\
$a$ & $0$ & $0$ & $0$ & $0$ & $0$ &$0$\\
$b$ & $0$ & $0$ & $b$ & $b$ & $b$ &$b$\\
$c$ & $0$ & $0$ & $b$ & $b$ & $b$ &$b$\\
$d$ & $0$ & $0$ & $b$ & $b$ & $b$ &$b$\\
$1$ & $0$ & $0$ & $b$ & $b$ & $b$ &$b$\\

\end{tabular}
\end{table}

\begin{remark}\label{remark100}
\emph{If $f$ is an $\mathfrak{f}$-mapping on a complete lattice $L$, then the binary operator $T$ in Theorem \ref{pro01} is a strong t-subnorm since $T(1,1)=f(1)\wedge f(1)=f(1)=1$.}
\end{remark}
\begin{theorem}\label{The01}
Let $L$ be a complete lattice whose top element $1$ is completely join-irreducible, and let $f$ be an $\mathfrak{f}$-mapping on $L$.
Then the following binary operator $T: L^2\rightarrow L$ is a left-continuous t-norm:
\begin{equation}\label{equation1}
T(x,y)=
\begin{cases}
x\wedge y &\emph{if } 1\in\{x,y\},\\
f(x)\wedge f(y) & \emph{otherwise}.
\end{cases}
\end{equation}
\end{theorem}
\proof
Clearly, $T$ is commutative with neutral element $1$. To verify that $T$ is non-decreasing, let $x, y, z\in L$ with $y\leq z$.
If $x=1$ then $T(x,y)=y\leq z=T(x,z)$. Suppose that $x\neq 1$. There are three cases as follows.

(i) If $y=1$ then $z=1$, which implies that $T(x,y)=x=T(x, z)$.

(ii) If $y\neq 1$ and $z=1$ then $T(x, y)=f(x)\wedge f(y)\leq x\wedge y\leq x=T(x,z)$
since $\mbox{Im}(f)$ is a subposet of $L$ and $f\leq \mbox{id}_L$.

(iii) If $y\neq 1$ and $z\neq 1$. Then $T(x, y)=f(x)\wedge f(y)\leq f(x)\wedge f(z)=T(x, z)$ since
$f$ is order-preserving.

Therefore, $T$ is non-decreasing.

Next, we show the associativity. Let $x, y, z\in L$. If $1\in \{x, y,z\}$ then $T(x, T(y,z))=T(T(x,y), z)$.
If $1\notin \{x,y,z\}$, then $f(y)\wedge f(z)\neq 1$. Thus
$$T(x,T(y,z))=T(x, f(y)\wedge f(z))=f(x)\wedge f(f(y)\wedge f(z))=f(x)\wedge f(y)\wedge f(z),$$
where the last equality follows from $f^2=f$ and $\mbox{Im}(f)$ forms a lattice. Moreover, by the commutativity of $T$, we also have that
$$T(T(x,y), z)=T(z, T(x, y))=f(z)\wedge f(x)\wedge f(y).$$
Hence, $T(x, T(y,z))=T(T(x,y), z)$, so $T$ is associative.

Finally, we only need to verify that $$T(x, \bigvee S)=\bigvee_{s\in S} T(x,s)$$ for all $x\in L$ and $S\subseteq L$ with $S\neq\emptyset$. If $x=1$ then
$T(x, \bigvee S)=\bigvee S=\bigvee_{s\in S} T(x,s)$. Assume that $x\neq 1$. Then there are two cases as follows.

Case 1. If $\bigvee S=1$ then $T(x, \bigvee S)=T(x, 1)=x$. Since $1$ is a completely join-irreducible element, $1\in S$.
On the other hand, we have
\begin{align*}
\bigvee_{s\in S} T(x,s)&=T(x, 1)\vee (\bigvee_{s\in S\setminus 1} T(x, s))\\
&=x\vee \bigvee_{s\in S\setminus 1} (f(x)\wedge f(s))\\
&=x\vee [f(x)\wedge (\bigvee_{s\in S\setminus 1} f(s))]\\
&\leq x\vee f(x)\\
&=x.
\end{align*}
Moreover, since $x\leq x\vee [f(x)\wedge \bigvee_{s\in S\setminus 1} f(s)]\leq x$, we conclude $\bigvee_{s\in S} T(x,s)=x$.
Hence, $$T(x, \bigvee S)=\bigvee_{s\in S} T(x,s).$$

Case 2. If $\bigvee S\neq 1$ then
\begin{align*}T(x, \bigvee S)&=f(x)\wedge f(\bigvee S)\\
&=f(x)\wedge (\bigvee_{s\in S} f(s))\\
&=\bigvee_{s\in S} (f(x)\wedge f(s))\\
&=\bigvee_{s\in S} T(x,s).
\end{align*}

To sum up, $T$ is a left-continuous t-norm.
\endproof

\begin{remark}\label{rema01}
\emph{In Theorem \ref{The01}, the condition ``the top element $1$ is completely join-irreducible" cannot be omitted.}
\end{remark}

\begin{example}\label{exa1}
\emph{Consider the finite lattice $L$ depicted in Fig. 1 again. Let the mapping $f: L\rightarrow L$ be defined as shown in Table \ref{Tab:01}.
It is verifiable that $f$ is an $\mathfrak{f}$-mapping. However, $1$ is not
a completely join-irreducible element. As is shown in Example \ref{fl}, there is no left-continuous t-norm on $L$.}
\end{example}
\begin{table}[!h]
\centering
\caption{A mapping $f$ on $L$}
\label{Tab:01}
\begin{tabular}{c|ccccccc}

 $x$ & $0$ & $a$ & $b$ & $c$ & $d$ & $e$ &$1$ \\
 \hline
$f(x)$ & $0$ & $0$ & $0$ & $c$ & $e$ & $e$ &$1$\\
\end{tabular}
\end{table}

\begin{example}\label{exa2}
\emph{Let $\mathbb{N}_0 \cup \{+\infty\}$ denote the linearly ordered lattice formed by extending the set of nonnegative
integers $\mathbb{N}_0$ with the positive infinity $+\infty$ in which the order relation is the natural extension of the standard order on
$\mathbb{N}_0$ and the lattice operations coincide with the usual maximum and minimum operations, respectively. Therefore, $+\infty$ serves as the top element of the lattice, which is not completely join-irreducible.
Define a mapping $f: \mathbb{N}_0 \cup \{+\infty\}\rightarrow \mathbb{N}_0 \cup \{+\infty\}$ as follows:}
\begin{equation*}
f(n)=
\begin{cases}
0 &\mbox{if } n\leq 1,\\
n& \mbox{if } +\infty>n\geq 2,\\
+\infty & \mbox{if } n=+\infty.
\end{cases}
\end{equation*}
\emph{Obviously, $f$ is an $\mathfrak{f}$-mapping on $\mathbb{N}_0 \cup \{+\infty\}$. Suppose that $T$ defined by \eqref{equation1} is a left-continuous t-norm on  $\mathbb{N}_0 \cup \{+\infty\}$. Then}

\begin{align*}1&=T(1, +\infty)\\
&=T(1, \bigvee_{+\infty>n\geq 2}n)\\
&=\bigvee_{+\infty>n\geq 2}T(1,n)\\
&=\bigvee_{+\infty>n\geq 2}(f(1)\wedge f(n))\\
&=\bigvee_{+\infty>n\geq 2}(0\wedge n)\\
&=0<1,
\end{align*}
\emph{a contradiction.}
\end{example}

From Remark \ref{rem00} and Theorem \ref{The01}, we obtain the following corollary.

\begin{corollary}\label{lemma01}
	Let $L$ be a complete lattice in which the top element $1$ is completely join-irreducible. Then there exists
	at least one $\mathfrak{f}$-mapping on $L$, which can induce a left-continuous t-norm
	on $L$.
\end{corollary}
\begin{proposition}
Let $f$ be an $\mathfrak{f}$-mapping on a complete lattice $L$ and $T: L^2\rightarrow L$ be the binary operator defined in Theorem \ref{The01}.
Then an element $x\in L$ is an idempotent element of $T$ if and only if $x\in \emph{Im}(f)$.
\end{proposition}
\proof
If $x\in \mbox{Im}(f)\setminus 1$, then there exists a $y\in L\setminus 1$ such that $f(y)=x$. Since $f$ is idempotent, we have $f(x)=f(f(y))=f(y)=x$. This implies
that $T(x, x)=f(x)\wedge f(x)=f(x)=x$. Note that $1=f(1)\in \mbox{Im}(f)$ and $T(1, 1)=1$. Therefore, for any $x\in \mbox{Im}(f)$, $T(x, x)=x$,
which means that $x$ is an idempotent element of $T$.

Conversely, suppose that $x$ is an idempotent element of $T$. Since $1\in \mbox{Im}(f)$, we only need to consider $x\in L\setminus 1$.
In this case, $x=T(x, x)=f(x)\wedge f(x)=f(x)\in \mbox{Im}(f)$, which implies that $x\in \mbox{Im}(f)$.
\endproof
\section{Constructing left-continuous t-norms}
When constructing left-continuous t-norms on complete lattices $L$, a natural approach is to construct
such t-norms via binary operators defined on certain subintervals of $L$. The main purpose of this section is to explore the 
construction of left-continuous t-norms on complete lattices in which every element is comparable with the endpoints of their 
respective specified subintervals.

Denote $\mathbb{N}=\{1,2,\cdots, n\}$ for any positive integer $n$.
Let $L$ be a complete lattice and $C=\{c_1, c_2, \cdots, c_n\}$ with $0=c_1< c_2< \cdots <c_n=1$ be a finite chain of $L$.
Then $L$ is a linear sum of $\{[c_i, c_{i+1}] \mid i\in \mathbb{N}\setminus n\}$ if and only if
$$L=\{x\in L\mid x \mbox{ and }c_i \mbox{ are comparable for all } i\in \mathbb{N}\}.$$

\begin{definition}\label{Def00001}
\emph{Let $L$ be a bounded lattice. Then a binary operator $T: L^2\rightarrow L$ is called annihilating if $T(0, y)= 0$
for all $y\in L$.}
\end{definition}

Clearly, every t-norm (resp. t-subnorm) on bounded lattices is annihilating. Specifically, the binary operator $\wedge$ on bounded lattices is also annihilating.

\begin{theorem}\label{the1.1}
Let a complete lattice $L$ be a linear sum of $\{[c_i, c_{i+1}] \mid i\in \mathbb{N}\setminus n\}$, and let $T_i$ be an annihilating
binary operator on $[c_i, c_{i+1}]$ for each $i\in \mathbb{N}\setminus n$.
Define a binary operator $T: L^2\rightarrow L$ by
\begin{equation}\label{equa1}
T(x,y)=
\begin{cases}
T_i(x,y) &\emph{if } (x,y)\in (c_i, c_{i+1}]^2,\\
x\wedge y& \emph{otherwise}.
\end{cases}
\end{equation}
Then $T$ is a left-continuous t-norm if and only if $T_i$ is a left-continuous t-subnorm for any $i\in \mathbb{N}\setminus \{n-1, n\}$, and $T_{n-1}$ is a left-continuous t-norm.
\end{theorem}
\proof
Suppose that $T$ is a left-continuous t-norm on $L$. By Eq.\eqref{equa1}, for all $i\in \mathbb{N}\setminus n$, we have
\begin{equation*}
T_i(x,y)=
\begin{cases}
T(x,y) &\mbox{if } (x,y)\in (c_i, c_{i+1}]^2,\\
c_i & \mbox{otherwise}
\end{cases}
\end{equation*}
since $T_i$ is annihilating on $[c_i, c_{i+1}]$. It is straightforward to verify that $T_i$ is a left-continuous t-subnorm for
all $i\in \mathbb{N}\setminus  n$. Furthermore, for any $x\in (c_{n-1}, c_{n}]$, $T_{n-1}(x,c_n)=T(x, c_n)=T(x,1)=x$ since $c_n=1$, which together
with $T_{n-1}$ being a left-continuous t-subnorm implies that $T_{n-1}$ is a left-continuous t-norm.

Conversely, suppose that $T_{n-1}$ is
a left-continuous t-norm on $[c_{n-1}, c_{n}]$, and $T_i$ is a left-continuous t-subnorm on $[c_i, c_{i+1}]$ for any $i\in \mathbb{N}\setminus \{n-1, n\}$ . Note that $c_n=1$.
Clearly, $T$ is commutative with neutral element $1$.

Firstly, we prove that $T$ is non-decreasing, i.e., for any $x,y,z\in L$, if $y\leq z$ then $T(x,y)\leq T(x,z)$. We consider the following five cases.

Case 1. If there exists an $i\in \mathbb{N}\setminus n$ such that $x,y,z\in (c_i, c_{i+1}]$, then $T(x,y)=T_i(x,y)\leq T_i(x,z)=T(x,z)$.

Case 2. If $x,y\in (c_i, c_{i+1}]$ for some $i\in \mathbb{N}\setminus n$, but $z\notin (c_i, c_{i+1}]$. Then $z>c_{i+1}\geq x$. This implies that
$T(x,y)=T_i(x,y)\leq x\wedge y\leq x= x\wedge z=T(x,z)$.

Case 3. If $x,z\in (c_i, c_{i+1}]$ for some $i\in \mathbb{N}\setminus n$, but $y\notin (c_i, c_{i+1}]$. Then $y\leq c_i< x$, which implies that
$T(x,y)=x\wedge y=y\leq c_i\leq T_i(x,z)=T(x,z)$.

Case 4. If $y,z\in (c_i, c_{i+1}]$ for some $i\in\mathbb{N}\setminus n$, but $x\notin (c_i, c_{i+1}]$. Then
$T(x,y)=x\wedge y\leq x\wedge z =T(x,z)$.

Case 5. If $x, y, z$ belong to three distinct intervals, respectively. Then $T(x,y)=x\wedge y\leq x\wedge z =T(x,z)$.

Thus $T$ is non-decreasing.

Secondly, we verify that $T$ is associative, i.e., $T(T(x,y), z)=T(x, T(y,z))$ for any $x, y, z\in L$. The proof is divided into three cases.

Case a. If $x,y,z\in (c_i, c_{i+1}]$ for some $i\in \mathbb{N}\setminus n$, then $T(T(x,y), z)=T(x, T(y,z))$ holds trivially.

Case b. Suppose that there are exactly two elements of $x,y,z$ which belong to $(c_i, c_{i+1}]$ for some $i\in \mathbb{N}\setminus n$, say $x$ and $y$. If $z> c_{i+1}$ then
\begin{align*}
T(T(x,y),z)&=T(T_i(x, y), z)\\
&=T_i(x, y)\wedge z\\
&=T_i(x,y)\\
&=T(x, y)\\
&=T(x, y\wedge z)\\
&=T(x, T(y,z)).
\end{align*}
If $z\leq c_i$ then
\begin{align*}
T(T(x,y),z)&=T(T_i(x, y), z)\\
&=T_i(x, y)\wedge z\\
&=z\\
&=T(x, z)\\
&=T(x, y\wedge z)\\
&=T(x, T(y,z)).
\end{align*}

Case c. If $x, y, z$ belong to three distinct intervals, respectively, then $$T(T(x,y),z)=T(x\wedge y,z)=x\wedge y\wedge z=T(x, y\wedge z)=T(x, T(y,z)).$$

Combining Cases a, b and c, the binary operator $T$ is associative.

Finally, we show that $T$ is left-continuous, i.e., for any $x\in L$ and $\emptyset\neq S\subseteq L$, $$T(x, \bigvee S)=\bigvee_{s\in S}T(x,s).$$
There are two cases as follows.

Case i. If $x, \bigvee S\in (c_i, c_{i+1}]$ for some $i\in \mathbb{N}\setminus n$, then we claim that there exists at least one element $y\in S$ such that
$y> c_i$. Otherwise, $\bigvee S\leq \bigvee [0, c_i]=c_i$, which implies that $\bigvee S \notin (c_i, c_{i+1}]$, a contradiction. Now, let $U=\{x\in S\mid x\leq c_i\}$.
Then $S\setminus U\neq\emptyset$ and $\bigvee S=\bigvee S\setminus U$. Clearly, $T(x, u)=x\wedge u\leq c_i$ for any $u\in U$ and
$T(x, v)=T_i(x, v)\geq c_i$ for any $v\in S\setminus U$. Thus
\begin{align*}
T(x,\bigvee S)&=T(x, \bigvee S\setminus U)\\
&=T_i(x, \bigvee S\setminus U)\\
&=\bigvee_{v\in S\setminus U} T_i(x,v)  \ \ (\mbox{by the left-continuity of } T_i)\\
&=[\bigvee_{u\in U}T(x, u)] \vee [\bigvee_{v\in S\setminus U} T_i(x,v)] \\
&=[\bigvee_{u\in U}T(x, u)] \vee [\bigvee_{v\in S\setminus U} T(x,v)] \\
&=\bigvee_{s\in S}T(x,s).
\end{align*}

Case ii. If $x,\bigvee S$ belong to distinct intervals, respectively, say $x\in (c_i, c_{i+1}], \bigvee S\in (c_j, c_{j+1}]$ with $i\neq j \in \mathbb{N}\setminus n$.

If $i< j$, then $x\leq  c_{i+1}\leq c_j < \bigvee S$. Similar to the proof of Case i, there exists an element $y\in S$ such that $x\leq  c_{i+1}\leq c_j < y$.
Note that $T(x, s)\leq x$ for any $s\in S$. Thus
\begin{align*}
T(x,\bigvee S)&=x\wedge (\bigvee S)\\
&=x\\
&\geq \bigvee_{s\in S} T(x,s)\\
&=[\bigvee_{z\in S\setminus y}T(x, z)] \vee T(x,y) \\
&=[\bigvee_{z\in S\setminus y}T(x, z)] \vee x \\
&\geq x,
\end{align*}
which means that $T(x,\bigvee S)=\bigvee_{s\in S}T(x,s)$.

If $i>j$, then  $x> c_{i}\geq c_{j+1} \geq \bigvee S\geq s$ for any $s\in S$. Thus
$$T(x, \bigvee S)=x\wedge (\bigvee S)=\bigvee S=\bigvee_{s\in S} (x\wedge s)=\bigvee_{s\in S}T(x,s).$$

Combining Cases i and ii, $T$ is left-continuous.

To sum up, the binary operator $T$ is a left-continuous t-norm.
\endproof

From Theorems \ref{pro01} and \ref{the1.1}, we have

\begin{corollary}\label{coro1}
Let a complete lattice $L$ be a linear sum of $\{[c_i, c_{i+1}] \mid i\in \mathbb{N}\setminus n\}$ and $T_{n-1}$ an annihilating binary operator on $[c_{n-1}, c_{n}]$.
Define a binary operator $T: L^2\rightarrow L$ by
\begin{equation*}
T(x,y)=
\begin{cases}
T_{n-1}(x,y) &\mbox{if } (x,y)\in (c_{n-1}, c_{n}]^2,\\
f_i(x)\wedge f_i(y) &\mbox{if } (x,y)\in (c_{i}, c_{i+1}]^2, 1\leq i\leq n-2,\\
x\wedge y & \mbox{otherwise},
\end{cases}
\end{equation*}
where every $f_i$ is a weak $\mathfrak{f}$-mapping on $[c_{i}, c_{i+1}]$ with $i\in \mathbb{N}\setminus \{n, n-1\}$.
Then $T$ is a left-continuous t-norm if and only if $T_{n-1}$ is a left-continuous t-norm.
\end{corollary}

\begin{definition}\label{Def01}
\emph{Let $\{[a_i, b_i]\mid i\in \mathbb{N}\}$ be a family of subintervals of a lattice $L$. We call $L$ a
semi-linear sum of $\{[a_i, b_i] \mid i\in \mathbb{N}\}$ if the following four conditions hold: for all $i,j\in \mathbb{N}$,\\
(i1) $(a_i,b_i]\neq \emptyset$;\\
(i2) $b_i\leq a_j$ if $i\leq j$ ;\\
(i3) $(a_i, b_i)\cap (a_j, b_j)=\emptyset$ if $i\neq j$;\\
(i4) $L=\{x\in L\mid \mbox{for all } i\in  \mathbb{N}, \mbox{ both } a_i \mbox{ and } b_i \mbox{ are comparable with } x\}$.}
\end{definition}

\begin{remark}\label{remark3.2}
\emph{Let $\{[a_i, b_i]\mid i\in \mathbb{N}\}$ be a family of subintervals of a complete lattice $L$.}

\emph{(1) Suppose that $a_1=0$ and $b_n=1$. If $i+1=j$ implies that $b_i=a_j$, then $L$ is a linear sum of $\{[a_i, b_i] \mid i\in \mathbb{N}\}$.}

\emph{(2) Let $T: L^2\rightarrow L$ be a left-continuous t-norm. If $T$ is closed on $(a_i, b_i]$ with $i\in \mathbb{N}$,
then $T|_{(a_i, b_i]}$ is a left-continuous t-subnorm.
In particular, if $T$ is closed on $(a_n, b_n]$ and $b_n=1$ then $T|_{(a_n, b_n]}$ is a left-continuous t-norm. If $T$ is closed on $(b_n, 1]$ and $b_n< 1$ then $T|_{(b_n, 1]}$ is a left-continuous t-norm.}

\emph{(3) Let $L$ be a semi-linear sum of $\{[a_i, b_i] \mid i\in \mathbb{N}\}$, and $x,y\in L$.
If $x$ and $y$ are non-comparable then both $x$ and $y$ belong to one of the following four type sets:
$$[a_i, b_i]\mbox{ with }i\in \mathbb{N}, [b_i, a_{i+1}]\mbox{ with }i\in \mathbb{N}\setminus n, [0, a_1]\mbox{ and }[b_n, 1].$$}
\end{remark}

Notice that by Remark \ref{remark3.2} (2), if a series of binary operators on subintervals are used to
construct left-continuous t-norms, then these operators must be left-continuous t-subnorms. In particular, the binary operator on the subinterval
containing the top element $1$ is required to be a left-continuous t-norm.

\begin{theorem}\label{the010}
Let a complete lattice $L$ be a semi-linear sum of $\{[a_i, b_i] \mid i\in \mathbb{N}\}$ and every $T_i$ be an 
annihilating binary operator on $[a_i, b_i]$ with $i\in \mathbb{N}$.
Define a binary operator $T: L^2\rightarrow L$ by
\begin{equation}\label{eqation2}
T(x,y)=
\begin{cases}
T_{n+1}(x,y) &\mbox{if } b_n<1 \mbox{ and } (x,y)\in (b_n, 1]^2,\\
T_i(x,y) &\mbox{if } (x,y)\in (a_i, b_i]^2,\\
f_i(x)\wedge f_i(y) &\mbox{if }(x,y)\in (b_{i}, a_{i+1}]^2, \\
f_0(x)\wedge f_0(y) &\mbox{if } (x,y)\in (0, a_{1}]^2,\\
x\wedge y& \mbox{otherwise},
\end{cases}
\end{equation}
where $f_0$ is a weak $\mathfrak{f}$-mapping on $[0, a_{1}]$, every $f_i$ is a weak $\mathfrak{f}$-mapping on $[b_i, a_{i+1}]$ with $i\in \mathbb{N}\setminus n$ and $T_{n+1}$ is an annihilating binary operator on $[b_n,1]$. Then the following three statements hold:
\begin{enumerate}[label=(\arabic*)]
\item If $b_n<1$, then $T$ is a left-continuous t-norm if and only if every $T_i$ is a left-continuous t-subnorm
with $i\in \mathbb{N}$ and $T_{n+1}$ is a left-continuous t-norm;
\item If $b_n=1$, then $T$ is a left-continuous t-norm if and only if every $T_i$ is a left-continuous t-subnorm with $i\in \mathbb{N}\setminus n$ and $T_{n}$ is a left-continuous t-norm.
\item If $T: L^2\rightarrow L$ is a left-continuous t-norm, then 
\begin{enumerate}[label=(\roman*)]
\item for all $i\in \mathbb{N}\setminus 1$, $T|_{[a_i, b_i]}=T_i$ if and only if 
$f_{i-1}$ is an $\mathfrak{f}$-mapping when $b_{i-1}< a_i$ and $T_{i-1}$ is a left-continuous strong t-subnorm when $b_{i-1}= a_i$.
\item if $0< a_1$, 
then $T|_{[a_1, b_1]}=T_1$ if and only if $f_{0}$ is an $\mathfrak{f}$-mapping.
\end{enumerate}
\end{enumerate}
\end{theorem}
\proof (1) and (2) are respective direct consequences of Theorems \ref{pro01} and \ref{the1.1}. Next, we just prove (3).
	
(i). Suppose that $T|_{[a_i, b_i]}=T_i$ for some $i\in \mathbb{N}\setminus 1$. Then $T(a_i, a_i)=T_i(a_i, a_i)=a_i$. 
	If $b_{i-1}< a_i$ then $T(a_i, a_i)=f_{i-1}(a_i)\wedge f_{i-1}(a_i)=f_{i-1}(a_i)$, which together with $T(a_i, a_i)=a_i$ implies $f_{i-1}(a_i)=a_i$. Thus $f_{i-1}$ is an $\mathfrak{f}$-mapping 
	when $b_{i-1}< a_i$. If $b_{i-1}= a_i$ then $T(a_i, a_i)=T(b_{i-1}, b_{i-1})=T_{i-1}(b_{i-1}, b_{i-1})$, which together with 
	$T(a_i, a_i)=a_i=b_{i-1}$ implies $T_{i-1}(b_{i-1}, b_{i-1})=b_{i-1}$. 
	Thus by (1) and (2), $T_{i-1}$ is a left-continuous strong t-subnorm when $b_{i-1}= a_i$. 
	
	Conversely, we only need to show that $T(x, y)=T_i(x,y)$ for all $(x, y)\in [a_i, b_i]^2$. The proof is divided into three cases.
	
	Case I. If both $x\neq a_i$ and $y\neq a_i$ then $(x,y)\in (a_i, b_i]^2$, which means that $T(x, y)=T_i(x, y)$. 
	
	Case II. If there is exactly one element of $x, y$ which equals to $a_i$, say $x=a_i$, then $y> a_i$. 
	Thus $$T(x, y)=T(a_i, y)=a_i\wedge y=a_i=T_i(a_i, y)=T_i(x,y)$$ since $T_i$ is annihilating.
	
	Case III. If $x=y=a_i$ then $T_i(x,y)=T_i(a_i, a_i)=a_i$. If $b_{i-1}< a_i$ then $$T(x, y)=T(a_i, a_i)=f_{i-1}(a_i)\wedge f_{i-1}(a_i)=a_i=T_i(x,y)$$ 
	since $f_{i-1}$ is an $\mathfrak{f}$-mapping. 
	If $b_{i-1}= a_i$ then $$T(x, y)=T(a_i, a_i)=T(b_{i-1}, b_{i-1})=T_{i-1}(b_{i-1}, b_{i-1})=b_{i-1}=a_i=T_i(x,y)$$ since $T_{i-1}$ is a strong t-subnorm.
	
	To sum up, $T|_{[a_i, b_i]}=T_i$.
	
	(ii). It is shown by complete analogy with (i).
	\endproof

The following two examples are used to illustrate Theorem \ref{the010}.
\begin{example}\label{example1}
\emph{Consider the complete lattice $L$ depicted in Fig. 3. Clearly, $L$ is a semi-linear sum of $\{[b, g]\}$.
Suppose that $T_1$ is a binary operator on $[b,g]$ as is shown in Table \ref{Tab:03} and $f_0$ is a mapping on $[0, b]$ as is shown in Table \ref{Tab:04}. It can be verified that $T_1$ is a left-continuous t-subnorm
and $f_0$ is a weak $\mathfrak{f}$-mapping. Let $T_2(x,y)=x\wedge y$ for all $x,y\in [g,1]$. Then
$T_2$ is a left-continuous t-norm on $[g,1]$. Then by Theorem \ref{the010}, $T: L^2\rightarrow L$ depicted in Table \ref{Tab:05} is a left-continuous t-norm.}
\end{example}

\par\noindent\vskip20pt
\begin{minipage}{10pc}
	\setlength{\unitlength}{0.7pt}\begin{picture}(600,330)
		\put(290,120){\circle{6}}\put(296,120){\makebox(0,0)[l]{\footnotesize $a$}}
		\put(290,80){\circle{6}}\put(286,68){\makebox(0,0)[l]{\footnotesize $0$}}
		\put(290,160){\circle{6}}\put(275,158){\makebox(0,0)[l]{\footnotesize $b$}}
		\put(290,240){\circle{6}}\put(302,240){\makebox(0,0)[l]{\footnotesize $g$}}
		\put(330,200){\circle{6}}\put(335,198){\makebox(0,0)[l]{\footnotesize $c$}}
		\put(250,200){\circle{6}}\put(235,198){\makebox(0,0)[l]{\footnotesize $d$}}
		\put(290,320){\circle{6}}\put(288,330){\makebox(0,0)[l]{\footnotesize $1$}}
		\put(330,280){\circle{6}}\put(335,278){\makebox(0,0)[l]{\footnotesize $e$}}
		\put(250,280){\circle{6}}\put(235,278){\makebox(0,0)[l]{\footnotesize $f$}}
		\put(290,200){\circle{6}}\put(275,198){\makebox(0,0)[l]{\footnotesize $h$}}
		\put(290,123){\line(0,1){34}}
		\put(290,163){\line(0,1){34}}
		\put(290,203){\line(0,1){34}}
		\put(292,162){\line(1,1){36}}
		\put(290,83){\line(0,1){34}}
		\put(292,242){\line(1,1){36}}
		\put(288,162){\line(-1,1){36}}
		\put(288,242){\line(-1,1){36}}
		\put(288,162){\line(-1,1){36}}
		\put(252,282){\line(1,1){36}}
		\put(328,282){\line(-1,1){36}}
		\put(252,202){\line(1,1){36}}
		\put(328,202){\line(-1,1){36}}
		\put(283,40){$L$}
		\put(170,20){ Fig.3. A finite lattice $L$.}
	\end{picture}
\end{minipage}

\begin{table}[!h]
\centering
\caption{A left-continuous t-subnorm $T_1$}
\label{Tab:03}
\begin{tabular}{c|ccccc}

 $T_1$ & $b$ & $d$ & $h$ & $c$ & $g$\\
 \hline
$b$ & $b$ & $b$ & $b$ & $b$ & $b$ \\
$d$ & $b$ & $b$ & $b$ & $b$ & $b$ \\
$h$ & $b$ & $b$ & $b$ & $b$ & $b$ \\
$c$ & $b$ & $b$ & $b$ & $b$ & $b$ \\
$g$ & $b$ & $b$ & $b$ & $b$ & $b$ \\
\end{tabular}
\end{table}

\begin{table}[!h]
\centering
\caption{A weak $\mathfrak{f}$-mapping $f_0$}
\label{Tab:04}
\begin{tabular}{c|ccccc}

 $x$ & $0$ & $a$ & $b$\\
 \hline
$f_0 $ & $0$ & $0$ & $a$ \\

\end{tabular}
\end{table}

\begin{table}[!h]
\centering
\caption{A left-continuous t-norm $T$}
\label{Tab:05}
\begin{tabular}{c|cccccccccc}

 $T$ & $0$ & $a$ & $b$ & $c$ & $d$& $h$ & $g$ & $e$ & $f$ & $1$\\
 \hline
 $0$ & $0$ & $0$ & $0$ & $0$ & $0$& $0$ & $0$ & $0$ & $0$ & $0$ \\
 $a$ & $0$ & $0$ & $a$ & $a$ & $a$& $a$ & $a$ & $a$ & $a$ & $a$ \\
 $b$ & $0$ & $a$ & $a$ & $b$ & $b$& $b$ & $b$ & $b$ & $b$ & $b$ \\
 $c$ & $0$ & $a$ & $b$ & $b$ & $b$& $b$ & $b$ & $c$ & $c$ & $c$ \\
 $d$ & $0$ & $a$ & $b$ & $b$ & $b$& $b$ & $b$ & $d$ & $d$ & $d$ \\
 $h$ & $0$ & $a$ & $b$ & $b$ & $b$& $b$ & $b$ & $h$ & $h$ & $h$\\
 $g$ & $0$ & $a$ & $b$ & $b$ & $b$& $b$ & $b$ & $g$ & $g$ & $g$\\
 $e$ & $0$ & $a$ & $b$ & $c$ & $d$& $h$ & $g$ & $e$ & $g$ & $e$\\
 $f$ & $0$ & $a$ & $b$ & $c$ & $d$& $h$ & $g$ & $g$ & $f$ & $f$\\
 $1$ & $0$ & $a$ & $b$ & $c$ & $d$& $h$ & $g$ & $e$ & $f$ & $1$\\
\end{tabular}
\end{table}

\begin{example}\label{example2}
\emph{Consider the complete lattice $L$ depicted in Fig. 3 again. Clearly, $L$ is
also a semi-linear sum of $\{[0, a], [g, 1]\}$. Suppose that both $T_1$ and $T_2$ are the strongest t-norm on $[0, a]$ and $[g, 1]$, respectively.
Let $f_1$ be a mapping on $[a, g]$ as is shown in table \ref{Tab:08}. It can be verified that $T_1$ is a left-continuous t-subnorm, $T_2$ is
a left-continuous t-norm and $f_1$ is a weak $\mathfrak{f}$-mapping. Then by Theorem \ref{the010}, $T: L^2\rightarrow L$ depicted in Table \ref{Tab:09} is a left-continuous t-norm.}
\end{example}

\begin{table}[!h]
\centering
\caption{A weak $\mathfrak{f}$-mapping $f_1$}
\label{Tab:08}
\begin{tabular}{c|cccccccc}

 $x$ & $a$ & $b$& $d$ & $h$ & $c$ & $g$\\
 \hline
$f_1$ & $a$ & $b$ & $b$ & $b$ & $b$ & $b$\\

\end{tabular}
\end{table}

\begin{table}[!h]
\centering
\caption{The left-continuous t-norm $T$}
\label{Tab:09}
\begin{tabular}{c|cccccccccc}

 $T$ & $0$ & $a$ & $b$ & $c$ & $d$& $h$ & $g$ & $e$ & $f$ & $1$\\
 \hline
 $0$ & $0$ & $0$ & $0$ & $0$ & $0$& $0$ & $0$ & $0$ & $0$ & $0$ \\
 $a$ & $0$ & $a$ & $a$ & $a$ & $a$& $a$ & $a$ & $a$ & $a$ & $a$ \\
 $b$ & $0$ & $a$ & $b$ & $b$ & $b$& $b$ & $b$ & $b$ & $b$ & $b$ \\
 $c$ & $0$ & $a$ & $b$ & $b$ & $b$& $b$ & $b$ & $c$ & $c$ & $c$ \\
 $d$ & $0$ & $a$ & $b$ & $b$ & $b$& $b$ & $b$ & $d$ & $d$ & $d$ \\
 $h$ & $0$ & $a$ & $b$ & $b$ & $b$& $b$ & $b$ & $h$ & $h$ & $h$\\
 $g$ & $0$ & $a$ & $b$ & $b$ & $b$& $b$ & $b$ & $g$ & $g$ & $g$\\
 $e$ & $0$ & $a$ & $b$ & $c$ & $d$& $h$ & $g$ & $e$ & $g$ & $e$\\
 $f$ & $0$ & $a$ & $b$ & $c$ & $d$& $h$ & $g$ & $g$ & $f$ & $f$\\
 $1$ & $0$ & $a$ & $b$ & $c$ & $d$& $h$ & $g$ & $e$ & $f$ & $1$\\
\end{tabular}
\end{table}

Note that it follows from Example \ref{example1} that $T(b,b)=a\neq b=T_1(b,b)$, which implies that $T|_{[b,g]}\neq T_1$.

\section{Conclusions}
This article focuses on the construction of left-continuous t-norms on complete lattices, and its core results are summarized as follows:

(i) Firstly, we introduced the concepts of $\mathfrak{f}$-mappings and weak $\mathfrak{f}$-mappings on complete lattices, respectively, and subsequently show that weak $f$-mappings 
induce left-continuous t-subnorms. Furthermore, $f$-mappings give rise to left-continuous t-norms 
on complete lattices whenever its top element $1$ is completely join-irreducible.

(ii) Secondly, we established some necessary and sufficient conditions for constructing left-continuous t-norms on complete lattices 
which are semi-linear sums of subinterval families via the ordinal sum of annihilating binary operators on each subinterval. 
Specifically, a left-continuous t-norm must be defined on the 
subinterval containing the top element $1$, whereas left-continuous t-subnorms are just required on the remaining subintervals.

It is regrettable that we do not delve further into the construction of left-continuous t-norms on complete lattices 
which do not admit a representation of semi-linear sums of the given subintervals.

\end{document}